\documentclass[preprintnumbers,amsmath,amssymb]{revtex4}

\usepackage{graphicx}
\usepackage{bm}

\newtheorem{theorem}{Theorem}[section]
\newtheorem{proposition}{Proposition}[section]
\newtheorem{lemma}{Lemma}[section]

\begin{document}

\title{Existence and congruence of global attractors for damped and forced integrable and nonintegrable discrete nonlinear Schr\"odinger equations}

\author{Dirk Hennig}

\date{\today}

\begin{abstract}
We study two damped and forced discrete nonlinear Schr\"odinger equations on the one-dimensional infinite lattice. 
Without damping and forcing they are represented by the integrable Ablowitz-Ladik equation (AL) featuring non-local cubic nonlinear terms, and its standard (nonintegrable) counterpart  with local cubic nonlinear terms (DNLS).   
The global existence of a unique solution to the initial value problem for both, the damped and forced AL and DNLS, is proven. It is further shown that for sufficiently close initial data, their corresponding solutions stay close for all times. Concerning the asymptotic behaviour of the solutions to the damped and forced AL and DNLS,  for the former a sufficient condition for the existence of a restricted global attractor is established while it is shown that the latter  possesses a global attractor. Finally, we prove the congruence of the restricted global  AL attractor and the  DNLS attractor 
for dynamics ensuing from  initial data contained in an appropriate bounded subset in a Banach space.
\end{abstract}

\maketitle
\section{Introduction}

The nonlinear Schr\"odinger equation (NLS) is one of the few examples of  completely integrable nonlinear partial differential equations.  Its study has lead to fundamental advances in 
 the theory and application of integrable systems possessing soliton solutions \cite{Ablowitz},\cite{Faddeev}. These soliton solutions  are of importance for a large number of physical and mathematical problems ranging from 
optical pulses propagation in nonlinear fibers to hydrodynamics, biophysics and condensed matter physics. 
 
However,  the application of continuum equations disregards often the inherent discrete lattice structure of the underlying system. The latter is naturally the case when dealing with sets of coupled units (oscillators) distributed e.g. in space,  and if,  the resulting phenomena of the nonlinear dynamics evolve on spatial scales comparable with the typical inter-unit distance. 

The issue of discretisation of NLS was addressed early in \cite{AL} where it was noticed that among a large number of possible discretisation, leading to nonlinear lattice systems, there is one that is also integrable termed the Ablowitz-Ladik, or AL equation. 
In contrast to the completely integrable AL equation, the so-called standard
discrete nonlinear Schr\"odinger equation, or DNLS equation, is known to be nonintegrable exhibiting also chaotic dynamics \cite{DNLS}.

Discrete nonlinear Schr\"odinger equations have been widely used in various contexts ranging from models of light propagation in arrays of optical waveguides \cite{application0}, the dynamics of atomic Bose-Einstein condensates trapped in optical lattices \cite{Morsch},  the study of denaturation of DNA double helix strands \cite{Peyrard}, breathers in granular materials \cite{Chong} to the dynamics of protein loops \cite{application1}. 

The AL equation, due to its complete integrability, has been utilised  for the development of perturbation theory exploring, for example, the problem of collisions between solitary waves in (nonintegrable) 
lattices \cite{Kevredikis1}. Notably, besides the soliton solutions another important class of solutions of the AL equation are rational solutions which are discrete versions of the Peregrine soliton and the Kuznetsov-Ma breather \cite{akhm_AL},\cite{akhm_AL2}. These rational solutions have become a field of intense research recently \cite{Kim2}-\cite{2DAL3}. 

For a more realistic description of many applications external forcing and dissipation needs to be included in the underlying model.   
Regarding the incorporation of dissipation (and external driving forces) into discrete systems, for the ensuing dissipative dynamics attractors for  infinite lattice dynamical systems have attracted considerable interest recently \cite{Bates}-\cite{Han}.
For their study the modern theory of infinite-dimensional dynamical systems provides powerful methods \cite{Hale}-\cite{Chueshov}.

We study the following two discrete nonlinear Schr\"odinger equations: 
\begin{equation}
i\frac{d \psi_n}{dt}=\kappa(\psi_{n+1}-2\psi_n+\psi_{n-1})+\mu\,|\psi_{n}|^2(\psi_{n+1}+\psi_{n-1})-i\delta \psi_n+g_n,\,\,\,n\in {\mathbb{Z}}\label{eq:AL}
\end{equation}
with $\psi_n \in {\mathbb{C}}$ and initial conditions:
\begin{equation}
 \psi_{n}(0)=\psi_{n,0},\,\,\,n \in {\mathbb{Z}},\label{eq:icsAL}
\end{equation}
and 
\begin{equation}
i\frac{d \phi_n}{dt}=\kappa(\phi_{n+1}-2\phi_n+\phi_{n-1})+ 
\gamma|\phi_n|^2\phi_n-i\delta \phi_n+g_n,\,\,\,n\in {\mathbb{Z}}\label{eq:DNLS}
\end{equation}
with $\phi_n \in {\mathbb{C}}$ and initial conditions:
\begin{equation}
 \phi_{n}(0)=\phi_{n,0},\,\,\,n \in {\mathbb{Z}}.\label{eq:icsDNLS}
\end{equation}
In what follows we refer to the damped and forced system (\ref{eq:AL}) and (\ref{eq:DNLS}) as dfAL and dfDNLS, respectively.
For the real parameter $\delta >0$, there is damping included in (\ref{eq:AL}),(\ref{eq:DNLS}), while $g=(g_n)_{n\in {\mathbb{Z}}}\neq 0$ serves in both systems as a general external force. The parameter $\mu$ and $\gamma$  determine the nonlinearity strength in (\ref{eq:AL}) and (\ref{eq:DNLS}), respectively.  The value of the parameter $\kappa \in \mathbb{R}$ 
regulates the coupling strength and without loss of generality we will set $\kappa=1$ subsequently.
In the conservative limit, for $\delta=0$, and without external force, i.e. $g=0$, the integrable AL equation and  and its nonintegrable DNLS counterpart result, respectively.
Notice that in (\ref{eq:AL}) and (\ref{eq:DNLS}) the nonlinear terms are both of cubic order. However, they are markedly different in the sense that, the nonlinear terms in (\ref{eq:AL}) are of nonlocal nature compared to the local terms in (\ref{eq:DNLS}). 
System (\ref{eq:DNLS}),{\ref{eq:icsDNLS}) 
with a general local nonlinear term has been studied in \cite{Nikos}. For more details 
concerning discrete nonlinear Schr\"odinger equations and their applications  we refer to \cite{DNLS},\cite{Kevrekidis}.

With the present work we study  an existence/closeness/congruence problem 
in the sense of  ``continuous dependence'' by investigating  closeness of the solutions of the dfAL and dfDNLS for close enough initial data. In particular the following questions are tackled: {\em (i) assuming that the initial data of the dfAL \eqref{eq:AL} and the dfDNLS \eqref{eq:DNLS} are sufficiently close in $l^2$, do the associated solutions remain close for sufficiently long times? 
(ii) While in the conservative limit the systems (\ref{eq:AL}) and (\ref{eq:DNLS}) exhibit distinct solution behaviour (integrability versus nonintegrability), the integrability of the AL equation gets destroyed by the inclusion of damping and external forcing. Then the question is, what do the two dissipative dynamics lattice  systems have in common? Does their asymptotic dynamics possess a global attractor and if so, what is the limit behaviour of the latter? Do the two system  even share a global attractor in the end?}
Moreover, from a wider point of view the answers to these questions seem  important because the closeness and congruence results  are not only relevant for discrete nonlinear Schr\"odinger equations but also for further investigations of the limit behaviour of  nonlinear (discrete) lattice systems in general.   For instance, the 
 asymptotic features of different discrete versions of physically important systems such as forced and damped continuum Ginzburg-Landau (GL) equations 
can be studied from the perspective described above (we refer to the section \ref{section:outlook} for more details). 

We answer the above questions by analytically proving that \textit{at least under certain smallness conditions on the initial data of the dfAL and the dfDNLS, the corresponding solutions remain close for sufficiently long times.} Crucially, with our congruence results we establish that the different discrete nonlinear Schr\"odinger equations exhibit the same asymptotic behaviour, that is they possess a  common global attractor.  
This is in as far interesting as  from the two discrete nonlinear Schr\"odinger equations without external forcing and damping  give rise to so profoundly different dynamics.

The outline of the  paper is as follows: We compare the solution properties of the two infinite lattice dynamical systems presented by the dfAL and its dfDNLS counterpart with particular attention to their asymptotic behaviour where we 
use methods from the theory of infinite-dimensional dynamical systems \cite{Hale}-\cite{Chueshov}. 
We establish the well-posedness by proving the global existence of a unique solution to the dfAL and dfDNLS. We demonstrate that, when the distance between the initial data for the dfAL and dfDNLS is sufficiently small in $l^2$ it remains small for sufficiently long time. 
We provide a sufficient criterion for the existence of a restricted global attractor for system (\ref{eq:AL}),(\ref{eq:icsAL}) and demonstrate the existence of a global attractor for system (\ref{eq:DNLS}),(\ref{eq:icsDNLS}). 
Finally, we prove the congruence of these two attractors, when the dynamics is initialised in an appropriate bounded subset in a Banach space. 


For the initial value problem for the infinite system of ordinary differential equations in (\ref{eq:AL}),(\ref{eq:icsAL}) and (\ref{eq:DNLS}),(\ref{eq:icsDNLS}), we consider solutions $\varphi=(\varphi_n)_{n\in {\mathbb{Z}}}\in C^1([0,\infty);l^2)$, where 
\begin{equation}
 l^2=\left\{ \varphi=(\varphi_n)_{n \in {\mathbb{Z}}}\,\in {\mathbb{C}}\,\,\,\vert \, \parallel \varphi\parallel_{l^2}=\left(\sum_{n \in {\mathbb{Z}}}|\varphi_n|^2\right)^{1/2}<\infty\right\}.\nonumber
\end{equation}

In the following  we use the notation
\begin{equation}
 B_R:=\left\{ \varphi \in l^2\,|\,
 \parallel \varphi\parallel_{l^2}<R\right\}\nonumber
\end{equation}
for the ball centered at $0$ of radius $R$  in $l^2$.

For any $\varphi \in l^2$ we define the linear operators $A,B,B^{*}:\,l^2 \rightarrow l^2$,
\begin{equation}
 (A\varphi)_{n }=\varphi_{n+1}-2\varphi_n+\varphi_{n-1},
\end{equation}
\begin{equation}
 (B\varphi)_{n}=\varphi_{n+1}-\varphi_n,\qquad (B^{*}\varphi)_{n}=\varphi_{n-1}-\varphi_n. \nonumber
\end{equation}
It holds that
\begin{equation}
 (B\varphi,\theta)_{l^2}=(\varphi,B^{*}\theta)_{l^2},\,\,\,\forall \varphi,\theta\in l^2,\nonumber
\end{equation}
and  $-A=BB^{*}=B^{*}B$, implying that 
\begin{equation}
 (A\varphi,\varphi)_{l^2}=-\parallel B\varphi\parallel_{l^2}\le 0,\,\,\,\forall \varphi \in l^2.\nonumber
\end{equation}
Furthermore, we observe that
\begin{equation}
 \parallel A\varphi\parallel_{l^2}\le 4\parallel \varphi\parallel_{l^2}.\nonumber
\end{equation}
As $A=A^*$, the  linear continuous operator $A$ is self-adjoint on $D(A)=l^2$ and $A\le 0$. Then the operator $A$ generates a uniformly continuous semigroup on $l^2$.

With the help of the transformation $\tilde{\varphi}(t)=\exp(-2it)\varphi(t)$ the linear operator $A$ will be replaced by
the linear operator $\Delta$ determined by  $(\Delta \varphi)_n=\varphi_{n+1}+\varphi_{n-1}$ in subsequent work. The effect of this transformation is  merely a shift of the continuous spectrum of $A$ so that it   comes to lie in the interval $[-2,2]$ instead of $[0,4]$ and one has 
\begin{equation}
 \parallel \Delta \varphi\parallel_{l^2}\le 2.\label{eq:boundDelta}
\end{equation} 
The advantage of this transformation is that with $\Delta$ a  more compact notation is achieved (see Eq.\,(\ref{eq:systemglobal})).

\section{Global existence of a unique solution}\label{subsection:existence}

For the current study of existence and uniqueness of a global solution of the dfAL and dfDNLS we combine them in a single system
\begin{equation}
i\frac{d \varphi_n}{dt}=(1+\mu\,|\varphi_{n}|^2)(\Delta \varphi)_n+ 
\gamma|\varphi_n|^2\varphi_n-i\delta \varphi_n(t)+g_n,\,\,\,n\in {\mathbb{Z}}\label{eq:systemglobal}
\end{equation}
with $\varphi_n \in {\mathbb{C}}$ and initial conditions:
\begin{equation}
 \varphi_{n}(0)=\varphi_{n,0},\,\,\,n \in {\mathbb{Z}}.\label{eq:icsglobal}
\end{equation}
Note that for $\gamma=0$ ($\mu=0$) the dfAL (dfDNLS) results from (\ref{eq:systemglobal}).

We formulate the infinite dimensional dynamical system (\ref{eq:systemglobal}),(\ref{eq:icsglobal}) as as an initial value problem in the Hilbert space $l^2$:
\begin{eqnarray}
 \dot{\varphi}&=&F(\varphi)\equiv -i(1+\mu\,|\varphi|^2)\Delta\varphi-i  
\gamma|\varphi|^2\varphi- \delta \varphi-i g,\,\,\,t>0,\label{eq:Hilbertsystem}\\
\varphi(0)&=&\varphi_0.\label{eq:Hilbertic}
\end{eqnarray}

We will use the following lemma in applying classical ODE theory:
\begin{lemma}
 \label{Lemma:Lipschitz}{\it \,\,Let $g=(g_n)_{n\in {\mathbb{Z}}}\in l^2$. The operator $F:\,l^2\rightarrow l^2$, defined by
\begin{equation}
\left(F(\theta)\right)_n= -i(1+\mu\,|\theta_{n}|^2)(\Delta \theta)_n-i
\gamma|\theta_n|^2\theta_n-\delta \theta_n(t)-ig_n\nonumber
\end{equation}
is Lipschitz continuous on bounded sets of $l^2$.}
\end{lemma}

\noindent{\bf Proof:} Let $\theta \in B_R$. 
For the nonlinear operator $N:l^2\rightarrow l^2$, $(N(\theta))_n=-i\mu\, |\theta_{n}|^2(\theta_{n+1}+\theta_{n-1})-i
\gamma|\theta_n|^2\theta_n$ we have 
\begin{eqnarray}
 \parallel N(\theta) \parallel_{l^2}^2&=&\sum_n \left|\mu\,|\theta_{n}|^2(\theta_{n+1}+\theta_{n-1})+
\gamma|\theta_n|^2\theta_n\right|^2\nonumber\\
&\le& \sum_n \left(\mu^2\,|\theta_{n}|^4|\theta_{n+1}+\theta_{n-1}|^2+
\gamma^2|\theta_n|^4|\theta_n|^2\right)\nonumber\\
&\le&(2\mu^2+\gamma^2)R^4\parallel \theta\parallel_{l^2}^2.\nonumber
\end{eqnarray}
Hence, $N:\,l^2\rightarrow l^2$ is bounded on bounded sets of $l^2$. 
For $\varphi,\theta \in B_R$ we derive 
\begin{eqnarray}
 \parallel N(\theta)-N(\varphi)\parallel_{l^2}^2&=&
 \sum_{n \in {\mathbb{Z}}}\left|\mu\,|\theta_{n}|^2(\theta_{n+1}+\theta_{n-1})+
\gamma|\theta_n|^2\theta_n-\mu\,|\varphi_{n}|^2(\varphi_{n+1}+\varphi_{n-1})-
\gamma|\varphi_n|^2\varphi_n\right|^2\nonumber\\
 &=&\frac{\mu^2}{4} \sum_{n \in {\mathbb{Z}}}| 
 \left(|\theta_n|^2+|\varphi_n|^2\right)\left[(\theta_{n+1}-\varphi_{n+1})+(\theta_{n-1}-\varphi_{n-1})\right]\nonumber\\
 &+&
 \left(|\theta_n|^2-|\varphi_n|^2\right)\left[(\theta_{n+1}+\varphi_{n+1})+(\theta_{n-1}+\varphi_{n-1})\right]|^2\nonumber\\
 &+&\frac{\gamma^2}{4}\sum_{n \in {\mathbb{Z}}}|
 (|\theta_n|^2 +|\varphi_n|^2)(\theta_n-\varphi_n)
 +(|\theta_n|^2 -|\varphi_n|^2)(\theta_n+\varphi_n)|^2
 \nonumber\\
 &\le& \left(\frac{\mu}{2}\right)^2\left\{
 2\sup_{n \in {\mathbb{Z}}}(|\theta_n|^2+|\varphi_n|^2)^2\,\sum_{n \in {\mathbb{Z}}}|\theta_n-\varphi_n|^2\right.\nonumber\\
 &+&\left.\sup_{n \in {\mathbb{Z}}}\left((|\theta_{n+1} +\varphi_{n+1}|^2+|\theta_{n-1} +\varphi_{n-1}|^2)
 \cdot \left(|\theta_n|+|\varphi_n|\right)^2 \right)\,\sum_{n \in {\mathbb{Z}}}|\theta_n-\varphi_n|^2
\right\}\nonumber\\
&+&\left(\frac{\gamma}{2}\right)^2 \left\{
\sup_{n \in {\mathbb{Z}}}(|\theta_n|^2 +|\varphi_n|^2)^2
\sum_{n \in {\mathbb{Z}}}
|\theta_n-\varphi_n|^2
 +\sup_{n \in {\mathbb{Z}}}(|\theta_n| +|\varphi_n|)^4
\sum_{n \in {\mathbb{Z}}}|\theta_n-\varphi_n|^2\right\}\nonumber\\
&\le& (6\mu^2+5\gamma^2)R^4\,\parallel \theta-\varphi\parallel^2_{l^2},\nonumber
\end{eqnarray}
verifying that the map $N:\,l^2\rightarrow l^2$ is Lipschitz continuous on bounded sets of $l^2$ with Lipschitz constant 
$L(R)=(6\mu^2+5\gamma^2)R^4$. Furthermore, since, due to (\ref{eq:boundDelta}), $\Delta$ is a bounded linear operator on $l^2$, we conclude that $F(\theta)$ is Lipschitz continuous on bounded sets of $l^2$.

\vspace*{0.5cm}

\hspace{16.5cm} $\square$

\vspace*{0.5cm}

For the proof of global existence of a unique solution to (\ref{eq:Hilbertsystem}),(\ref{eq:Hilbertic}) for $\mu \neq 0, \gamma=0$ we will use the following statement:


\begin{lemma}
 {\it \,\,If $P=\sum_{n\in{\mathbb{Z}}}\ln(1+|\psi_n|^2)< \infty$, then $\parallel \psi\parallel_{l^2}<\infty$.}
\label{Lemma:Pmu}
\end{lemma}

\noindent{\bf Proof:} We have 
\begin{equation}
  P=\sum_{n\in{\mathbb{Z}}}\ln(1+|\psi_n|^2)< \infty\,\,\,\Leftrightarrow\,\,\, \forall \epsilon >0\,\,\,\exists N_{\epsilon} \in \mathbb{N}\,\,\,s.t.\,\,\,\sum_{|n|\ge N_{\epsilon}}\ln(1+|\psi_n|^2)< \epsilon.\nonumber
\end{equation}

Writing $\lambda_n=\ln(1+|\psi_n|^2)$ 
and choosing $N_{\epsilon}$ such that $\epsilon <1/2$ we get 
\begin{eqnarray}
 \sum_{|n|\ge N_{\epsilon}}|\psi_n|^2&=&\sum_{|n|\ge N_{\epsilon}}\left(\exp(|\lambda_n|^2)-1\right)=\sum_{|n|\ge N_{\epsilon}}\,\left(\sum_{k=0}^{\infty}\frac{|\lambda_n|^{2k}}{k!}-1\right)\nonumber\\
 &=&\sum_{|n|\ge N_{\epsilon}}\,\sum_{k=1}^{\infty}\frac{|\lambda_n|^{2k}}{k!}=\sum_{k=1}^{\infty}\sum_{|n|\ge N_{\epsilon}}\,\frac{|\lambda_n|^{2k}}{k!}\nonumber\\
&\le& 
\sum_{k=1}^{\infty}\sum_{|n|\ge N_{\epsilon}}\,|\lambda_n|^{2k}
\le \sum_{k=1}^{\infty}\left(\sum_{|n|\ge N_{\epsilon}}\,|\lambda_n|^{2}\right)^k\nonumber\\
 &\le&\sum_{k=1}^{\infty}\epsilon^{k}=\frac{1}{1-\epsilon}-1< 2\epsilon .\nonumber
\end{eqnarray}
In conclusion, for all $0\le \epsilon<1/2$ there exists $0<N_{\epsilon}\in {\mathbb{N}}$ such that 
\begin{equation}
 \sum_{|n|\ge N_{\epsilon}}|\psi_n|^2<2\epsilon,\,\,\,\forall t\ge 0,\nonumber
\end{equation}
and hence, $\sum_{n\in{\mathbb{Z}}}|\psi_n|^2<\infty$.

\vspace*{0.5cm}

\hspace{16.5cm} $\square$

\vspace{0.5cm}

Regarding the global existence of a unique solution to (\ref{eq:Hilbertsystem}),(\ref{eq:Hilbertic})  we have the following:

\begin{proposition}
 {\it \,\,For every $\varphi_0\in l^2$ 
the system (\ref{eq:Hilbertsystem}) possesses a unique global solution $\varphi(t)$ on $[0,\infty)$ belonging to
$C^1([0,\infty),l^2)$.}\label{Proposition:unique}
\end{proposition}

\noindent{\bf Proof:} 
With the proven Lipschitz continuity of the operator $F$ on bounded sets of $l^2$, for any given initial data $\varphi_0\in l^2$  
the existence of a unique solution $\varphi(t)\in C^1([0,T_0),l^2)$ for some $T_0>0$ can be verified by standard methods from the theory of ODEs (see e.g. \cite{Zeidler}). 
Whenever $T_0<\infty$ then 
$\lim_{t\rightarrow T_0^{-}}\parallel \varphi(t)\parallel_{l^2}=\infty$. 

\vspace{0.5cm}
Moreover, we show the global existence of the solutions, that is, $T_0=\infty$. 
The ensuing analysis is appropriately performed separately for the dfAL and dfDNLS.  
First we treat the dfAL for which for convenience we set $\mu=1$. (Note that by the transformation ${\tilde\varphi}(t)=\sqrt{\mu} \varphi(t)$ the amplitude can be accordingly rescaled.)  
For the dfAL we consider 
\begin{eqnarray}
 \frac{d}{dt}{P}(t)&=&\frac{d}{dt}\sum_{n\in{\mathbb{Z}}}\ln(1+|\psi_n(t)|^2)= 2\sum_{n\in{\mathbb{Z}}}\frac{
 {\rm Re} \psi_n(t) {\rm Im} g_n -{\rm Re} g_n {\rm Im} \psi_n(t) -\delta |\psi_n(t)|^2}{ 1+|\psi_n(t)|^2}\nonumber\\
 &\le& -2\delta \sum_{n\in{\mathbb{Z}}} \frac{|\psi_n(t)|^2}{1+|\psi_n(t)|^2}+4\sum_{n\in{\mathbb{Z}}} \frac{|g_n||\psi_n(t)|}{ 1+|\psi_n(t)|^2}\nonumber\\
 &\le& -2\delta \sum_{n\in{\mathbb{Z}}} \frac{|\psi_n(t)|^2}{1+|\psi_n(t)|^2}+
 4\sum_{n\in{\mathbb{Z}}} \left[ \frac{\delta}{2}\frac{|\psi_n(t)|^2}{ (1+|\psi_n(t)|^2)^2}+\frac{1}{2\delta}|g_n|^2\right]\nonumber\\
 &=&-2\delta \sum_{n\in{\mathbb{Z}}}\left( \frac{|\psi_n(t)|^2}{ 1+|\psi_n(t)|^2}\right)^2+\frac{2}{\delta}\parallel g\parallel_{l^2}^2 < \infty,\,\,\,\forall t\ge 0,\nonumber
\end{eqnarray}
where we used Young's inequality.
Therefore, $P(t)< \infty$ for all $t\ge 0$.
Hence, by  Lemma \ref{Lemma:Pmu} above we obtain that $\parallel \psi(t)\parallel_{l^2}
< \infty$ for all $t\ge 0$.

To demonstrate global existence for the dfDNLS we consider 
\begin{eqnarray}
 \frac{d}{dt} \parallel\phi(t)\parallel_{l^2}^2  &=&2\sum_{n\in{\mathbb{Z}}}[{\rm Re} \phi_n(t) {\rm Im} g_n-
 {\rm Re} g_n {\rm Im} \phi_n(t)]
 -\delta |\phi_n(t)|^2\nonumber\\
&\le& -2\delta \sum_{n\in{\mathbb{Z}}} |\phi_n(t)|^2+4\sum_{n\in{\mathbb{Z}}}|g_n||\phi_n(t)|\nonumber\\
&\le& -\delta \sum_{n\in{\mathbb{Z}}} |\phi_n(t)|^2+ \frac{4}{\delta}\parallel g\parallel_{l^2}^2 < \infty,\,\,\,\forall t\ge 0.\nonumber
\end{eqnarray}
That is,  $\parallel \phi(t)\parallel_{l^2}^2
< \infty$ for all $t\ge 0$.

\vspace*{0.5cm}

\hspace{16.5cm} $\square$

\vspace*{0.5cm}
In conclusion, for the dfAL for any $\psi_0 \in l^2$, as well as for the dfDNLS for any 
$\phi_0 \in l^2$,  
the corresponding solution $\psi (t)$  of (\ref{eq:AL}),(\ref{eq:icsAL}) and $\phi (t)$ of  (\ref{eq:DNLS}),(\ref{eq:icsDNLS}), respectively, 
is bounded 
for all $t \in [0,\infty)$.  The solution operator determined by 
\begin{equation}
 S_{\mu}(t):\psi_0 \in l^2 \rightarrow  \psi(t)=S_{\mu}(t)\psi_0\in l^2,\,\,\,t\ge 0,\nonumber
\end{equation}
and 
\begin{equation}
 S_{\gamma}(t):\phi_0 \in l^2 \rightarrow  \phi(t)=S_{\gamma}(t)\phi_0\in l^2,\,\,\,t\ge 0,\nonumber
\end{equation}
generates a continuous semigroup $\left\{ S_{\mu}(t)\right\}_{t\ge 0}$ and 
 $\left\{ S_{\gamma}(t)\right\}_{t\ge 0}$ on $l^2$, respectively.

\section{Closeness of the dfAL and dfDNLS solutions}

In this section we demonstrate that for sufficiently small initial conditions for the dfAL and its dfDNLS counterpart  (\ref{eq:icsAL}) and (\ref{eq:icsDNLS}), respectively, and provided the $l^2-$distance between them is sufficiently small, 
the distance between the associated solutions to (\ref{eq:AL}) and (\ref{eq:DNLS}) remains small for sufficiently long  times $t>0$.

\begin{theorem}
{\it \,\,For every $t_0>0$ there exist a small  $\epsilon_0>0$ and positive constants $C_0$ and $C$ such that for every $\epsilon \in (0,\epsilon_0)$ for which  the initial conditions of the dfAL, $\psi(0)$,   and the dfDNLS, $\phi(0)$, and $g\in l^2$ satisfy:
\begin{equation}
\parallel \phi(0)-\psi(0)\parallel_{l^2}\le C_0 \epsilon^3,\label{eq:distance0}
\end{equation}
and
\begin{equation}
\frac{2}{\delta}\parallel g\parallel_{l^2}\le \parallel \phi(0)\parallel_{l^2}\le K_{\phi}\epsilon,\label{eq:Cphi}
\end{equation}
with some constant $K_{\phi}>0$, the corresponding solutions fulfill for every $t\in [0,t_0]$ 
\begin{equation}
 \parallel \phi(t)-\psi(t)\parallel_{l^2}\le C \epsilon^3.\label{eq:boundy}
\end{equation}
}\label{Theorem:closeness}
\end{theorem}

\noindent{\bf Proof:} 
Introducing for the (local) distance the variable 
$y_n=\psi_n-\phi_n$ one derives using the Cauchy-Schwarz inequality
\begin{eqnarray}
  \frac{d}{dt}\parallel y(t)\parallel_{l^2}^2&=&2\parallel y(t)\parallel_{l^2} \frac{d}{dt}\parallel y(t)\parallel_{l^2}\nonumber\\
  &=& \sum_{n \in {\mathbb{Z}}}\left\{i[(\overline{y}_{n+1}+
   \overline{y}_{n-1})y_n-({y}_{n+1}+
   {y}_{n-1})\overline{y}_n]-2\delta |y_n|^2\right.\nonumber\\   
   &+&\left.i\mu |y_n+\phi_n|^2[(\overline{y}_{n+1}+\overline{\phi}_{n+1}
   +\overline{y}_{n-1}+\overline{\phi}_{n-1})y_n-
   ({y}_{n+1}+{\phi}_{n+1}
   +{y}_{n-1}+{\phi}_{n-1})\overline{y}_n]\right.\nonumber\\
   &-&\left.2i\gamma |\phi_n|^2(\overline{\phi}_n y_n-\phi_n \overline{y}_n)\right\}\nonumber\\   
   &\le& (4-2\delta)\parallel y(t)\parallel_{l^2}^2+
   4\mu\left(\parallel y(t)\parallel_{l^2}^3+\parallel \phi(t)\parallel_{l^2}^2 \parallel y(t)\parallel_{l^2}
   +\parallel y(t)\parallel_{l^2}^2 \parallel \phi(t)\parallel_{l^2}+\parallel \phi(t)\parallel_{l^2}^3\right) \parallel y(t)\parallel_{l^2}
\nonumber\\
   &+&2\gamma \parallel \phi(t)\parallel_{l^2}^3 \parallel y(t)\parallel_{l^2}\nonumber
\end{eqnarray}
where  we exploited the continuous embeddings
 $l^r\subset l^s,\,\,\,\parallel \phi\parallel_{l^s}\le \parallel\phi\parallel_{l^r},\,\,\,1 \le r\le s \le \infty$.
Let $t_0>0$ be given. Define 
\begin{equation}
 \overline{T}_0 =\sup \left\{\overline{t} \in [0,t_0]: \sup_{t\in [0,\overline{t}]}\parallel y(t)\parallel_{l^2} \le K_y\epsilon^3\right\}.\label{eq:timeinterval}
\end{equation}
In Lemma \ref{Lemma:asymboundgamma} below  we establish that if $\frac{2}{\delta}\parallel g\parallel_{l^2}\le \parallel \phi(0)\parallel_{l^2}$ one has 
$\parallel \phi(t)\parallel_{l^2}\le \parallel \phi(0)\parallel_{l^2}$ for all $t>0$.

Then we obtain
 \begin{eqnarray}
  \frac{d}{dt}\parallel y(t)\parallel_{l^2}&\le& 
  (2-\delta)\parallel y(t)\parallel_{l^2}\nonumber\\
  &+&\left(2\mu \left(K_y^3\epsilon^6+K_{\phi}^2 K_y\epsilon^2+K_y^2 K_{\phi} \epsilon^4+K_{\phi}^3\right)+\gamma K_{\phi}^3\right)\epsilon^3.\nonumber
\end{eqnarray}
Furthermore, for every $t\in [0,\overline{T}_0]$ and sufficiently small $\epsilon >0$ one can find a positive constant $M_1$  independent of $\epsilon$ such that
\begin{equation}
2\mu \left(K_y^3\epsilon^6+K_{\phi}^2 K_y\epsilon^2+K_y^2 K_{\phi} \epsilon^4+K_{\phi}^3\right)+\gamma K_{\phi}^3\le M_1,\nonumber
\end{equation}
giving
\begin{equation}
 \frac{d}{dt}\parallel y(t)\parallel_{l^2}\le \epsilon^3 M_1+M_2\parallel y(t)\parallel_{l^2},\nonumber
\end{equation}
with $M_2=2-\delta$.
Applying Gronwall's inequality one gets
\begin{equation}
 \frac{d}{dt}\parallel y(t)\parallel_{l^2}\exp(-M_2 t)-\parallel y(0) \parallel_{l^2}\le \int_0^t \epsilon^3 M_1 \exp(-M_2 s) ds \le \epsilon^3 \frac{M_1}{M_2},\nonumber
\end{equation}
yielding with the assumption $\parallel y(0)\parallel_{l^2}\le C_0 \epsilon^3$ for every 
$t\in [0,\overline{T}_0]$:
\begin{equation}
 \parallel y(t)\parallel_{l^2}\le \left(C_0+\frac{M_1}{M_2}\right)\exp(M_2 t)\epsilon^3.\nonumber
\end{equation}
Thus one can set
\begin{equation}
 K_y=\left(C_0+\frac{M_1}{M_2}\right)\exp(M_2 t),\nonumber
\end{equation}
and the time interval in  (\ref{eq:timeinterval}) can be extended to the entire time range with $\overline{T}_0=t_0$ using an elementary continuation argument. 
This concludes the proof.

\vspace*{0.5cm}

\hspace{16.5cm} $\square$

\vspace*{0.5cm}

After having shown that the distance between the solutions of the 
dfAL and the dfDNLS measured in terms of the $l^2-$metric remains small (bounded above by ${\cal{O}}(\epsilon^3)$), compared to the $l^2-$norm of the solutions themselves (bounded above by ${\cal{O}}(\epsilon$), we demonstrate analogous features for the $l^{\infty}-$norm 
determining the maximal distance between individual units.

\begin{theorem}
{\it \,\,Let the assumption of Theorem \ref{Theorem:closeness} hold. Assume further that 
\begin{equation}
\parallel y(0)\parallel_{l^2} \le\parallel y(0)\parallel_{l^1}\le L_1\epsilon^3\label{eq:yl1}
\end{equation}

For every $t_0>0$, there exist a small  $\epsilon_0>0$ and a positive constant $C_{\infty}$  such that for every $\epsilon \in (0,\epsilon_0)$ for 
the maximal distance satisfies
\begin{equation}
 \parallel \phi(t)-\psi(t)\parallel_{l^\infty}\le C_{\infty}\, \epsilon^3,\,\,\,t\in[0,t_0].\label{eq:boundy3}
\end{equation}
}\end{theorem}

\noindent{\bf Proof:} For the time evolution of the distance variable $y_n=\psi_n-\phi_n$  we derive
\begin{eqnarray}
  \frac{d}{dt}y_n&=&
   -i({y}_{n+1}+{y}_{n-1})-\delta y_n-i\left[\mu|y_n+\phi_n|^2({y}_{n+1}+{y}_{n-1}+{\phi}_{n+1}+{\phi}_{n-1})-\gamma |\phi_n|^2\phi_n\right].\nonumber
\end{eqnarray}   

Performing a spatial Fourier transform 
\begin{equation}
 y_n(t)=\frac{1}{2\pi}\int_{-\pi}^{\pi}\hat{y}_q(t)\exp(i qn) dq,\qquad\hat{y}_q(t)=\sum_n y_n(t)\exp(-i q n),\nonumber
\end{equation}
gives the system
\begin{equation}
 \dot{\hat{y}}_q(t)=-(2i\cos q+\delta)\hat{y}_q(t)-i\hat{F}_q(t),\nonumber
\end{equation}
the formal solution of which is given by 
\begin{equation}
 \hat{y}_q(t)=\hat{y}_q(0)\exp[-(2i\cos q+\delta)t]-i\int_0^t \hat{F}_q(\tau) \exp[-(2i\cos q+\delta)(t-\tau)]d\tau.\label{eq:FT}
\end{equation}
Applying the inverse Fourier transform to (\ref{eq:FT}) 
we obtain
\begin{equation}
 y_n(t)=\frac{1}{2\pi}\int_{-\pi}^{\pi}\hat{y}_q(0)\exp[-(2i\cos q+\delta)t]\exp(iqn)dq
 -\frac{i}{2\pi}\int_{-\pi}^{\pi}\int_0^t \hat{F}_q(\tau) \exp[-(2i\cos q+\delta)(t-\tau)]d\tau \exp(iqn)dq.\nonumber
\end{equation}
With the help of the Cauchy-Schwarz inequality and $\parallel y\parallel_{l^4}\le \parallel y\parallel_{l^2}$ 
we estimate as follows:
\begin{eqnarray}
  |y_n(t)|&\le&\frac{1}{2\pi} \left|\int_{-\pi}^{\pi} \hat{y}_q(0) \exp[-(2i\cos q+\delta)t] \exp(iqn) dq \right|\nonumber\\
 &+&  \frac{1}{2\pi} \left|\int_{-\pi}^{\pi}\int_0^t \hat{F}_q(\tau) \exp[-(2i\cos q+\delta)(t-\tau)] d\tau \exp(iqn)dq \right|\nonumber\\
 &\le& \frac{1}{2\pi} \int_{-\pi}^{\pi} \left|\hat{y}_q(0)\right| dq
 +\frac{1}{2\pi} \int_{-\pi}^{\pi} \int_0^t \left|\hat{F}_q(\tau)\right|d\tau dq\nonumber\\
  &=&\frac{1}{2\pi}\int_{-\pi}^{\pi} \left|\sum_n y_n(0)\exp(-iqn)\right| dq
 +\frac{1}{2\pi}\int_{-\pi}^{\pi}\int_0^t \left|\sum_n F_n(\tau)\exp(-iqn)\right|d\tau dq\nonumber\\
  &\le& 
  \frac{1}{2\pi}\int_{-\pi}^{\pi} \sum_n \left| y_n(0)\right| dq+
  \frac{1}{2\pi}\int_{-\pi}^{\pi}\int_0^t \sum_n \left| F_n(\tau)\right|d\tau dq=\sum_n \left| y_n(0)\right| +
   \int_0^t \sum_n \left| F_n(\tau)\right|d\tau\nonumber\\
  &=&\parallel y(0)\parallel_{l^1}+ 
  \int_0^t\sum_n  \left|
  \mu|y_n(\tau)+\phi_n(\tau)|^2({y}_{n+1}(\tau)+{y}_{n-1}(\tau)+{\phi}_{n+1}(\tau)+{\phi}_{n-1}(\tau))\right.\nonumber\\
  &-&\left.\gamma |\phi_n(\tau)|^2\phi_n(\tau)\right|d\tau\nonumber\\
&\le&
  \parallel y(0)\parallel_{l^1}+ 
  \int_0^t\left[2\mu (\parallel y(\tau)\parallel_{l^4}^2 +\parallel \phi(\tau)\parallel_{l^4}^2) 
 (\parallel y(\tau)\parallel_{l^2}+\parallel \phi(\tau)\parallel_{l^2})+\gamma  \parallel \phi(\tau)\parallel_{l^4}^2 \parallel \phi(\tau)\parallel_{l^2}\right]\nonumber\\
 &\le&
  \parallel y(0)\parallel_{l^1}+ 
  \int_0^t\left[2\mu (\parallel y(\tau)\parallel_{l^2}^2 +\parallel \phi(\tau)\parallel_{l^2}^2) 
 (\parallel y(\tau)\parallel_{l^2}+\parallel \phi(\tau)\parallel_{l^2})+\gamma  \parallel \phi(\tau)\parallel_{l^2}^2 \parallel \phi(\tau)\parallel_{l^2}\right]\nonumber
 \end{eqnarray}

 Let $t_0>0$ be given and define 
\begin{equation}
 \overline{T}_0 =\sup \left\{\overline{t} \in [0,t_0]: \sup_{t\in [0,\overline{t}]}=\parallel y(t)\parallel_{l^2} \le K_y\epsilon^3\right\}.\label{eq:timeinterval1}
\end{equation}
Then we derive the following upper bound:
\begin{eqnarray}
|y_n(t)|&\le&\parallel y(0)\parallel_{l^1}+\left[2\mu(K_y^2\epsilon^4+K_{\phi}^2)(K_y\epsilon^2+K_{\phi})+\gamma K_{\phi}^3\right]t\cdot \epsilon^3,\qquad \forall t\in [0,\overline{T}_0].\nonumber
 \end{eqnarray}

For every $t\in [0,\overline{T}_0]$  and sufficiently small $\epsilon >0$ one can find a positive constant $M_1$  independent of $\epsilon$ such that
\begin{equation}
2\mu \left(K_y^3\epsilon^6+K_{\phi}^2 K_y\epsilon^2+K_y^2 K_{\phi} \epsilon^4+K_{\phi}^3\right)+\gamma K_{\phi}^3\le M_1,\nonumber
\end{equation}
giving with the hypothesis (\ref{eq:yl1}) for every 
$t\in [0,\overline{T}_0]$:
\begin{equation}
 |y_n(t)|\le (L_1+M_1)t\cdot\epsilon^3.\nonumber
\end{equation} 
Taking the supremum over $n \in {\mathbb{Z}}$ one gets
 \begin{equation}
  \sup_{n \in {\mathbb{Z}}}
  |y_n(t)|=\parallel y(t)\parallel_{l^{\infty}}\le (L_1+M_1)t\cdot\epsilon^3,\qquad \forall t\in [0,\overline{T}_0].\nonumber
 \end{equation}
Thus one can set
\begin{equation}
 K_{y,\infty}=(L_1+M_1)t,\nonumber
\end{equation}
and the time interval  
can be extended to the entire time range with $\overline{T}_0=t_0$ by an elementary continuation argument.

\vspace*{0.5cm}

\hspace{16.5cm} $\square$

We remark that it is desirable to obtain a growth rate of the 
 distance $||y(t)||_{l^2}$ that is uniformly bounded for any $\epsilon>0$ and any finite $t$ in $(0,\infty)$ as
\begin{eqnarray}
\label{gr}
\frac{d}{dt}||y(t)||_{l^2}\leq M\,\varepsilon^3,\nonumber
\end{eqnarray}
where $M$ depends on the parameters and initial data but not on $t$, and consequently, the distance between solutions grows at most linearly
for any $t\in (0,\infty)$, as
\begin{eqnarray}
\label{gr1}
||y(t)||_{l^2}\leq M\, t\,\varepsilon^3.\nonumber
\end{eqnarray}

\section{Existence of a restricted global attractor for the semigroup   $\left\{S_{\mu}(t)\right\}_{t\ge 0}$ in $l^2$}

 In order to distinguish in the following between  a global attractor (existent for the dfDNLS) and a restricted global attractor (relevant for the dfAL),  we recall their  definitions:  
 
 \vspace{0.5cm} 
 
 \noindent{\bf Definition:}{\it \,\,A set ${\cal A}_{\gamma}\subset l^2$ is called a global attractor for the semigroup  $\left\{S_{\gamma}(t)\right\}_{t\ge 0}$ associated with the system (\ref{eq:DNLS}),(\ref{eq:icsDNLS}) in $l^2$ if 
 (i)  ${\cal A}_{\gamma}\neq 0$ is a compact subset of $l^2$, (ii) an invariant set, that is,  
 $S_{\gamma}(t){\cal A}_{\gamma}={\cal A}_{\gamma}$ for all $t\ge 0$, 
  and, (iii) an attracting set for $\left\{S_{\gamma}\right\}_{t\ge 0}$ in $l^2$, that is, for all bounded $B\subset l^2$, it holds that 
 $\lim_{t\rightarrow \infty}\,dist(S_{\gamma}(t)B,{\cal A}_{\gamma})=0$, where the Hausdorff semi-distance between two nonempty subsets $U,V$ of $l^2$ is determined by
 \begin{equation}
  dist(U,V)=\sup_{u\in U}\,\inf_{v\in V}\,d(u,v)_{l^2}.\nonumber
 \end{equation}}

 \noindent{\bf Definition:}{\it \,\,Let $\left\{S_{\mu}(t)\right\}_{t\ge 0}$ be the semigroup associated with the system (\ref{eq:AL}),(\ref{eq:icsAL}). We say that ${\cal A}_{\mu}\subset l^2$ is a restricted global attractor for  $\left\{S_{\mu}(t)\right\}_{t\ge 0}$
 in $l^2$ if for some closed, nonempty subset $U\subset B_{R_{\mu}}$ of $l^2$, $S_{\mu}(t):U\rightarrow U$ ($t\ge 0$) is a  semigroup on $U$ such that ${\cal A}_{\mu}$ is a global attractor for $\left\{S_{\mu}(t)\right\}_{t\ge 0}$
restricted to $U$, that is, (i) $S_{\mu}(t){\cal A}_{\mu}={\cal A}_{\mu}$ for $t\ge 0$, (ii) ${\cal A}_{\mu}$ is compact, and (iii) ${\cal A}_{\mu}$ attracts solutions of bounded subsets of $U$ \cite{Hale}.}

\vspace{0.5cm}

In this section we establish the existence of a restricted global attractor for the semigroup  $\left\{S_{\mu}(t)\right\}_{t\ge 0}$ and in the next section the existence of a global attractor for  $\left\{S_{\gamma}(t)\right\}_{t\ge 0}$ is shown.

\vspace{0.5cm}

\subsection{Existence of an absorbing set in $l^2$}

 First we explore the existence of an absorbing set in $l^2$ for the dynamical system associated with the damped and forced dfAL (\ref{eq:AL}),(\ref{eq:icsAL}) in the asymptotic regime $t \rightarrow \infty$.
 
We need  the following lemma:
\begin{lemma}
 {\it \,\,Assume that $\delta$ and $g=(g_n)_{n\in {\mathbb{Z}}}$
 satisfy
 \begin{equation}
  \delta^2 < 3 \sum_{n\in {\mathbb{Z}}} |g_n|^{4/3}.\label{eq:conditionglobal}
 \end{equation}
 Then for any $\psi_0$ the solutions exist globally in time and are uniformly bounded satisfying
 \begin{equation}
  \parallel \psi(t) \parallel_{l^2}\le \parallel \psi(0) \parallel_{l^2}.\label{eq:uniformbound}
 \end{equation}
 }
 \end{lemma}

\noindent{\bf Proof:} 
From Lemma \ref{Lemma:Pmu} we get 
\begin{eqnarray}
 \frac{d}{dt}{P}(t)&=&
 \frac{d}{dt}\sum_{n\in{\mathbb{Z}}}\ln(1+|\psi_n(t)|^2)\nonumber\\
  &\le& -2\delta \sum_{n\in{\mathbb{Z}}} \frac{|\psi_n(t)|^2}{1+|\psi_n(t)|^2}+4\sum_{n\in{\mathbb{Z}}} \frac{|g_n||\psi_n(t)|}{ 1+|\psi_n(t)|^2}.\nonumber
 \end{eqnarray}
The relation (\ref{eq:uniformbound}) is satisfied if  $\dot{P}(t)\le 0$, which is true if  
 \begin{equation}
  2\sum_{n\in{\mathbb{Z}}} \left(\frac{1}{ 1+|\psi_n(t)|^2}(-\delta |\psi_n|^2+2|\psi_n||g_n|)\right)\le 0,\nonumber
 \end{equation}
 leading to
 \begin{equation}
  -2\delta \sum_{n\in{\mathbb{Z}}} |\psi_n|^2+4\sum_{n\in{\mathbb{Z}}}|\psi_n||g_n|\le 0.\label{eq:suff}
 \end{equation}
 With the aid of Young's inequality and the continuous embedding $l^2\subset \l^4$ one derives
 \begin{equation}
  \sum_{n\in{\mathbb{Z}}}|\psi_n||g_n|\le \frac{3}{4}\sum_{n\in{\mathbb{Z}}}|g_n|^{4/3}+\frac{1}{4}\sum_{n\in{\mathbb{Z}}}|\psi_n|^4\le \frac{3}{4}\sum_{n\in{\mathbb{Z}}}|g_n|^{4/3}+\frac{1}{4}\parallel \psi\parallel_{l^2}^4,\nonumber
 \end{equation}
so that one  obtains  the following sufficient condition 
 \begin{equation}
  -2\delta \parallel \psi\parallel_{l^2}^2+\parallel \psi\parallel_{l^2}^4+3\sum_{n\in{\mathbb{Z}}}|g_n|^{4/3}\le 0,\label{eq:ineqglobal}
 \end{equation}
for that the relation (\ref{eq:suff}) holds.

If condition (\ref{eq:conditionglobal}) holds, then the inequality (\ref{eq:ineqglobal}) is satisfied for any $\parallel \psi\parallel_{l^2}^2 \in {\mathbb{R}_+}$, so that  it follows that $\dot{P}(t)\le 0$, which yields the uniform bound (\ref{eq:uniformbound}}).

\vspace*{0.5cm}

\hspace{16.5cm} $\square$

Now we state the main assertion:

\begin{lemma}
 {\it \,\,Assume that the hypothesis (\ref{eq:conditionglobal}) holds and that 
\begin{equation}
 0<\parallel \psi(0)\parallel_{l^2}^2\le R_{\mu}^2<\frac{\delta}{4\mu}.\label{eq:ALbound}
\end{equation}
Let $(g_n)_{n\in {\mathbb{Z}}}=g \in l^2$ and $(\psi_n(0))_{n\in {\mathbb{Z}}}=\psi_0$. For the dynamical system determined by (\ref{eq:AL}),(\ref{eq:icsAL})
\begin{equation}
 S_{\mu}(t):\,\psi_0\in l^2\rightarrow \psi(t)\in l^2,
\end{equation}
there exists  a bounded absorbing set $B_{r}$ in $l^2$, that is, for every  set $B\subset B_{R_\mu}$ of $l^2$, there is a number $t_0(B,B_{r})>0$ such that  
$S_{\mu}(t)B\subset B_{r}$ for all $t\ge t_0(B,B_{r})$.}
\label{Lemma:absorbingmu}
\end{lemma}

\noindent{\bf Proof:} With the assumption (\ref{eq:ALbound}) we estimate
 \begin{eqnarray}
  \frac{d}{dt}\sum_{n\in {\mathbb{Z}}}|\psi_n|^2&=& 
  \sum_{n\in {\mathbb{Z}}}\left( 
  i\,\mu |\psi_n|^2 [\left(\overline{\psi}_{n+1}+\overline{\psi}_{n-1}\right)\psi_{n}- 
  \left({\psi}_{n+1}+{\psi}_{n-1}\right)\overline{\psi}_{n}
  ] +i\,(\overline{g}_n\psi_n-g_n\overline{\psi}_n)-2\delta \sum_{n\in {\mathbb{Z}}}|\psi_n|^2 \right)\nonumber\\
  &\le&4\mu \sup_{t\ge 0}\sup_{n\in{\mathbb{Z}}}|\psi_n(t)|^2
  \parallel \psi(t)\parallel_{l^2}^2-2\delta \parallel \psi(t)\parallel_{l^2}^2+4\parallel g\parallel_{l^2}\parallel \psi(t)\parallel_{l^2}\nonumber\\
  &\le& 4\mu \sup_{t\ge 0}\parallel\psi(t)\parallel_{l^\infty}^2
  \parallel \psi(t)\parallel_{l^2}^2-2\delta \parallel \psi(t)\parallel_{l^2}^2+4\parallel g\parallel_{l^2}\parallel \psi(t)\parallel_{l^2}\nonumber\\
  &\le& 4\mu \parallel\psi(0)\parallel_{l^2}^2
  \parallel \psi(t)\parallel_{l^2}^2-2\delta \parallel \psi(t)\parallel_{l^2}^2+4\parallel g\parallel_{l^2}\parallel \psi(t)\parallel_{l^2}\nonumber\\
  &\le& 4\mu R_{\mu}^2 \parallel \psi(t)\parallel_{l^2}^2-2\delta \parallel \psi(t)\parallel_{l^2}^2+4\parallel g\parallel_{l^2}\parallel \psi(t)\parallel_{l^2}.\label{eq:estimateattractorpsi}
 \end{eqnarray}
 Using Young's inequality for the last term on the right side of (\ref{eq:estimateattractorpsi})
 \begin{equation}
  4\parallel g\parallel_{l^2}\parallel \psi(t)\parallel_{l^2}\le 
  \frac{4}{\delta}\parallel g\parallel_{l^2}^2
 + \delta \parallel \psi(t)\parallel_{l^2}^2,
 \end{equation}
we have
\begin{equation}
  \frac{d}{dt}\parallel \psi(t)\parallel_{l^2}^2 +2\delta \parallel \psi(t)\parallel_{l^2}^2\le (4\mu R_{\mu}^2 +\delta)\parallel \psi(t)\parallel_{l^2}^2+\frac{4}{\delta}\parallel g\parallel_{l^2}^2,\nonumber
\end{equation}
so that 
\begin{equation}
  \frac{d}{dt}\parallel \psi(t)\parallel_{l^2}^2 +\left(\delta-4\mu R_{\mu}^2\right) \parallel \psi(t)\parallel_{l^2}^2\le \frac{4}{\delta}\parallel g\parallel_{l^2}^2,\label{eq:psiestimate}
\end{equation}
ensuring  that for $0<R_{\mu}^2<\delta/(4\mu)$ that $\psi \in L^{\infty}([0,\infty),l^2)$.
Applying Gronwall's inequality gives:
\begin{equation}
 \parallel \psi(t)\parallel_{l^2}^2\le \parallel \psi(0)\parallel_{l^2}^2\exp(-(\delta-4\mu R_{\mu}^2)t)+\frac{4}{\delta}\frac{\parallel g\parallel_{l^2}^2}{\delta-4\mu R_{\mu}^2}\left(1-\exp(-(\delta-4\mu R_{\mu}^2)t)\right).\label{eq:psiGronwall}
\end{equation}

In the asymptotic regime $t\rightarrow \infty$ this leads to
\begin{equation}
 {\lim\sup}_{t \rightarrow \infty} \parallel \psi(t)\parallel_{l^2}^2\le \frac{4}{\delta}\frac{\parallel g\parallel_{l^2}^2}{\delta-4\mu R_{\mu}^2}.\nonumber
\end{equation}
Defining 
\begin{equation}
 \rho^2=\frac{4}{\delta}\frac{\parallel g\parallel_{l^2}^2}{\delta-4\mu R_{\mu}^2},\nonumber
\end{equation}
we observe that for any number $r$, satisfying $R_{\mu}>r>\rho$, the ball $B_{r}$ of $l^2$ is an absorbing set for the semigroup $S_{\mu}(t):$ That is, for a  set  
$B\in B_R$  
it follows that, for $t\ge t_0(B,B_{r})$, where 
\begin{equation}
 t_0=\frac{1}{\delta-4\mu R^2}\log\left(\frac{R^2-\rho^2}{{r}^2-\rho^2}\right),\nonumber
\end{equation}
 one has $\parallel \psi(t)\parallel_{l^2}^2\le r^2$, i.e. 
$S_{\mu}(t)B\subset B_{r}$.

\vspace*{0.5cm}

\hspace{16.5cm} $\square$

\vspace{0.5cm}

Note that, although we ensured in section \ref{subsection:existence} the global existence of a unique solution to the dfAL in $l^2$, the nonlocal feature of the nonlinear term of the dfAL
allows to show the existence of an absorbing set only when the sufficient condition (\ref{eq:ALbound}) is satisfied with the effect that all solutions are contained 
in a ball $R_{\mu}$ in $l^2$ for all $t\ge 0$ (cf. Eq. (\ref{eq:estimateattractorpsi})).

\subsection{Asymptotic compactness of the semigroup $\left\{S_{\mu}(t)\right\}_{t\ge 0}$}

Here we verify that the semigroup $\left\{S_{\mu}(t)\right\}_{t\ge 0}$ associated with the dfAL (\ref{eq:AL}),(\ref{eq:icsAL}) possesses the asymptotic tail end property. 

\begin{lemma}
{\it \,\, } Let $(\psi_n(0))_{n\in {\mathbb{Z}}}=\psi_0 \in B$
and $(g_n)_{n\in {\mathbb{Z}}}=g \in l^2$. For any $\xi>0$ there exist
 $T(\xi)$ and $M(\xi)$ such that the solution $\psi(t)$ of (\ref{eq:AL}),(\ref{eq:icsAL}) satisfies for all $t\ge T(\xi)$:
 \begin{equation}
  \sum_{|n|> 2K}|\psi_n(t)|^2\le \xi\,\,\,\,{\rm for\, any\,\,\,\,} K>M(\xi).\label{eq:asymptotic}
 \end{equation}
\label{Lemma:asymtailmu}
\end{lemma}

\noindent{\bf Proof:} For a contradiction let us suppose that this assertion is   not true, 
i.e. there is an $\epsilon_0>0$ and a subsequence $(n_k)_{k \in {\mathbb{Z}}}$ of $\mathbb{Z}$ such that
\begin{equation}
 \sum_{|k|>m}|\psi_{n_k}(t)|^2 \ge \epsilon_0,\,\,\,\forall t\ge 0,\,\,\,{\rm for\,\, any}\,\, m\in {\mathbb{N}}.\nonumber
\end{equation}
Using (\ref{eq:psiestimate}) we have 
\begin{equation}
  \frac{d}{dt}\sum_{|k|>m}|\psi_{n_k}(t)|^2 +\left[\delta-4\mu R_{\mu}^2\right] \sum_{|k|>m}|\psi_{n_k}(t)|^2\le \frac{4}{\delta}  \sum_{|k|>m}|g_{n_k}|^2,\nonumber
\end{equation}
and Gronwall's inequality gives
\begin{equation}
 \sum_{|k|>m}|\psi_{n_k}(t)|^2 \le \sum_{|k|>m}|\psi_{n_k}(0)|^2 \exp(-(\delta-4\mu R_{\mu}^2)t)+\frac{4}{\delta}\frac{\sum_{|k|>m} |g_{n_k}|^2}{\delta-4\mu R_{\mu}^2}\left(1-\exp(-(\delta-4\mu R_{\mu}^2)t)\right).\label{eq:partialsum}
 \end{equation}
 From (\ref{eq:partialsum}) we infer
 \begin{equation}
  \inf_{t\in[0,\infty]}\sum_{|k|>m}|\psi_{n_k}(t)|^2 =
  \frac{4}{\delta}\frac{\sum_{|k|>m} |g_{n_k}|^2}{\delta-4\mu R^2} \ge \epsilon_0>0,\,\,\,{\rm for\,\, any}\,\, m\in {\mathbb{N}}.\nonumber
 \end{equation}
 Therefore, for every $m\in {\mathbb{N}}$, $\sum_{|k|>m}|g_{n_k}|^2>C(\delta,\mu,R_{\mu})\cdot \epsilon_0>0$,  which contradicts the fact that for $(g_n)_{n\in {\mathbb{Z}}}=g$ with $\parallel g  \parallel_{l^2}<\infty$, every subsequence $(g_{n_k})_{k\in {\mathbb{N}}}$ must  converge. 

\vspace*{0.5cm}

\hspace{16.5cm} $\square$

\vspace*{0.5cm} 

\noindent{\bf Definition:} The semigroup  $\left\{S_{\mu}(t)\right\}_{t\ge 0}$ is said to be asymptotically compact in $l^2$ if, for any bounded $B\subset B_{R_{\mu}} \subset  l^2$, and any sequence $\left\{t_n\right\}$,  $\left\{\phi_n\right\}$ with $t_n\ge 0$, $t_n \rightarrow \infty$ as $n\rightarrow \infty$, and $\phi_n \in B\subset B_{R_{\mu}}$, the sequence  $\left\{S_{\mu}(t_n)\phi_n\right\} $ is relatively compact in $l^2$.

\vspace{0.5cm}

\begin{proposition}
{\it \,\,Under the same conditions of  Lemma \ref{Lemma:asymtailmu} the semigroup $\left\{S_{\mu}(t)\right\}_{t\ge 0}$  is asymptotically compact.}
\label{Proposition:asymcompmu}
\end{proposition}

\noindent{\bf Proof:} 
By contradiction: Suppose that every subsequence of $\psi^n(t_n)=S_{\mu}(t_n)\psi_0^n \in B$ diverges in $l^2$ as $t_n \rightarrow \infty$ (equivalently, no sequence $\psi^n(t_n)=S_{\mu}(t_n)\psi_0^n$ has a convergent subsequence in $l^2$ as $t_n \rightarrow \infty$).
Take any subsequence  $\psi^{n_k}(t_{n_k})=S_{\mu}(t_{n_k})\psi_0^{n_k} \in B$, $k\in {\mathbb{N}}$. Suppose $\parallel  \psi^{n_k}(t_{n_k}) \parallel_{l^2}^2 \rightarrow \infty$ as $t_{n_k}\rightarrow \infty$, that is, $k\rightarrow \infty$. 

Then for any $M\in {\mathbb{R}}>0$ there are only finitely many $k$ such that 
\begin{equation}
 \parallel  \psi^{n_k}(t_{n_k}) \parallel_{l^2}^2\le M.\label{eq:Mk}
\end{equation}

 Denote all values of $k$ for which (\ref{eq:Mk}) is satisfied by $k_1,...,k_m$. Setting $N_M=\max\left\{k_1,...,k_m\right\}+1$, then for any $k>N_M$ it holds that
 \begin{equation}
  \parallel  \psi^{n_k}(t_{n_k}) \parallel_{l^2}^2=\sum_{i\in {\mathbb{Z}}}|\psi^{n_k}_i(t_{n_{k}})|^2>M,\,\,\,\forall k>N_M.\label{eq:normpsiM}
 \end{equation}
We  split the infinite sum in (\ref{eq:normpsiM}) as
\begin{equation}
\sum_{i\in {\mathbb{Z}}}|\psi^{n_k}_i(t_{n_{k}})|^2=\sum_{|i|>2L}|\psi^{n_k}_i(t_{n_{k}})|^2+\sum_{|i|\le 2L}|\psi^{n_k}_i(t_{n_{k}})|^2
>M,\,\,\,\forall k>N_M,\label{eq:normpsiM1}
 \end{equation}
 for any fixed $0<L<\infty$. For the finite sum we get $\sum_{|i|\le 2L}|\psi^{n_k}_i(t_{n_{k}})|^2<M_L$ with $M>M_L$, implying that
 \begin{equation}
\sum_{|i|\ge 2L}|\psi^{n_k}_i(t_{n_{k}})|^2
>M-M_L>0,\,\,\,\forall k>N_M.\label{eq:normpsiM2}
 \end{equation}
 Since the relation (\ref{eq:normpsiM2}) holds for all 
 $t\ge t_{n_{k> N_M}}$, for all $M>0$ and any $L>0$, it contradicts the asymptotic tail end property of $\left\{S_{\mu}(t)\right\}_{t\ge 0}$  as established by Lemma \ref{Lemma:asymtailmu}.

\vspace*{0.5cm}

\hspace{16.5cm} $\square$

Finally, facilitating Proposition \ref{Proposition:asymcompmu} and Theorem 1.1 in \cite{Temam}, we are now able to state the main result of this section.

\begin{theorem}
 {\it \,\,
 The semigroup $\left\{S_{\mu}(t)\right\}_{t\ge 0}$ associated with the dfAL  (\ref{eq:AL}),(\ref{eq:icsAL}) possesses a unique restricted global attractor ${\cal A}_{\mu}\subset B_{r}\subset l^2$.}\label{Theorem:attractor AL}
 \end{theorem}

\vspace{0.5cm}

\section{Existence of a global attractor for the semigroup $\left\{S_{\gamma}(t)\right\}_{t\ge 0}$}

Here we recall some results regarding the existence of an absorbing set in $l^2$ for the dynamical system belonging to the dfDNLS (\ref{eq:DNLS}),(\ref{eq:icsDNLS}) in the asymptotic regime $t \rightarrow \infty$ (see also \cite{Nikos}).

\begin{lemma}
{\it \,\,
Let $(g_n)_{n\in {\mathbb{Z}}}=g \in l^2$ and $(\phi_n(0))_{n\in {\mathbb{Z}}}=\phi_0$. For the dynamical system determined by (\ref{eq:DNLS}),(\ref{eq:icsDNLS})
\begin{equation}
 S_{\gamma}(t):\,\phi_0\in l^2\rightarrow \phi(t)\in l^2
\end{equation}
there exists  a bounded absorbing set $B_{\tilde{r}}$ in $l^2$, that is, for every bounded set $B$ of $l^2$, there is a $t_0(B,B_{\tilde{r}})$ such that  
$S_{\gamma}(t)B\subset B_{\tilde{r}}$ for all $t\ge t_0(B,B_{\tilde{r}})$. 
Furthermore, if
\begin{equation}
 \frac{2}{\delta}\parallel g\parallel_{l^2}\le \parallel \phi_0\parallel_{l^2},\label{eq:conddeltaphi0}
\end{equation}
then it holds 
\begin{equation}
 \parallel \phi(t)\parallel_{l^2}\le \parallel \phi_0\parallel_{l^2},\,\,\,\forall t\ge 0.\label{eq:phitbelow}
\end{equation}
\label{Lemma:asymboundgamma}}
\end{lemma}

\noindent{\bf Proof:} We derive for the dfDNLS 
\begin{eqnarray}
  \frac{d}{dt}\sum_{n\in {\mathbb{Z}}}|\phi_n|^2&=& 
  \sum_{n\in {\mathbb{Z}}}\left( i\,(\overline{g}_n\phi_n-g_n\overline{\phi}_n)-2\delta \sum_{n\in {\mathbb{Z}}}|\phi_n|^2 \right)\nonumber\\
  &\le& -2\delta \parallel \phi(t)\parallel_{l^2}^2+4\parallel g\parallel_{l^2}\parallel \phi(t)\parallel_{l^2}.\label{eq:estimateattractorphi}
 \end{eqnarray}
 This gives
\begin{equation}
  \frac{d}{dt}\parallel \phi(t)\parallel_{l^2}^2 +\delta \parallel \phi(t)\parallel_{l^2}^2\le \frac{4}{\delta}\parallel g\parallel_{l^2}^2,\nonumber
\end{equation}
from which with the use of Gronwall's inequality we obtain:
\begin{equation}
 \parallel \phi(t)\parallel_{l^2}^2\le \parallel \phi(0)\parallel_{l^2}^2\exp(-\delta t)+\frac{4}{\delta^2}\parallel g\parallel_{l^2}^2\left(1-\exp(-\delta t)\right).\label{eq:phit}
\end{equation}

Asymptotically for $t\rightarrow \infty$ it results that
\begin{equation}
 {\lim\sup}_{t \rightarrow \infty} \parallel \phi(t)\parallel_{l^2}^2\le \frac{4}{\delta^2}\parallel g\parallel_{l^2}^2.\nonumber
\end{equation}
Defining 
\begin{equation}
 \tilde{\rho}^2=\frac{4}{\delta^2}\parallel g\parallel_{l^2}^2,\nonumber
\end{equation}
we observe that for any number $\tilde{r}>\tilde{\rho}$ the ball $B_{\tilde{r}}\subset l^2$  
is an absorbing set for the semigroup $S_{\gamma}(t).$ That is, if $B$ is a bounded set of $l^2$ included in a ball $B_R$ 
it follows that for $t\ge t_0(B,B_{\tilde{r}})$ where 
\begin{equation}
 t_0=\frac{1}{\delta}\log\left(\frac{R^2-\tilde{\rho}^2}{\tilde{r}^2-\tilde{\rho}^2}\right),\nonumber
\end{equation}
one has $\parallel \psi(t)\parallel_{l^2}^2\le \tilde{r}^2$, that is, 
$S_{\gamma}(t)B\subset B_{\tilde{r}}$.

Finally, from the relations 
\begin{equation}
 \parallel \phi(t)\parallel_{l^2}^2\le \parallel \phi(0)\parallel_{l^2}^2\exp(-\delta t)+\frac{4}{\delta^2}\parallel g\parallel_{l^2}^2\left(1-\exp(-\delta t)\right)\le \parallel \phi(0)\parallel_{l^2}^2\nonumber
\end{equation}
one obtains (\ref{eq:phitbelow}).

\vspace*{0.5cm}

\hspace{16.5cm} $\square$

\vspace*{0.5cm}

Concerning the asymptotic tail end property, for  
the dfDNLS (\ref{eq:DNLS}),(\ref{eq:icsDNLS}) we have the following lemma. 

\vspace*{0.5cm}
\begin{lemma}
 {\it \,\,  Let $(\phi_n(0))_{n\in {\mathbb{Z}}}=\psi_0 \in B$, where $B$ is a bounded set of $l^2$ and 
 $(g_n)_{n\in {\mathbb{Z}}}=g \in l^2$. For any $\xi>0$ there exist
 $T(\xi)$ and $M(\xi)$ such that the solution $\phi(t)$ of (\ref{eq:DNLS}),(\ref{eq:icsDNLS}) satisfies for all $t\ge T(\xi)$:
 \begin{equation}
  \sum_{|n|\ge 2K}|\phi_n(t)|^2\le \xi\,\,\,\,{\rm for\, any\,\,\,\,} K>M(\xi).\label{eq:asymptoticDNLS}
 \end{equation}}\label{Lemma:asymtailgamma}
 \end{lemma}

\noindent{\bf Definition:} The semigroup  $\left\{S_{\gamma}(t)\right\}_{t\ge 0}$ is said to be asymptotically compact in $l^2$ if, for any bounded $B  \subset  l^2$, and any sequence $\left\{t_n\right\}$,  $\left\{\phi_n\right\}$ with $t_n\ge 0$, $t_n \rightarrow \infty$ as $n\rightarrow \infty$, and $\phi_n \in B$, the sequence  $\left\{S_{\gamma}(t_n)\phi_n\right\} $ is relatively compact in $l^2$.

\begin{proposition}
{\it \,\,Under the same conditions of Lemma \ref{Lemma:asymtailgamma},  the semigroup $\left\{S_{\gamma}(t)\right\}_{t\ge 0}$ is asymptotically compact.}
\label{Proposition:asymcompgamma}
\end{proposition}

The proofs of Lemma \ref{Lemma:asymtailgamma} and  Proposition \ref{Proposition:asymcompgamma}  proceed in a similar manner to the corresponding proofs for the dfAL and are omitted here.

In conclusion, by virtue of  Proposition \ref{Proposition:asymcompmu} and Theorem 1.1 in \cite{Temam}, we have:

\begin{theorem}
{\it \,\,
The semigroup $\left\{S_{\gamma}(t)\right\}_{t\ge 0}$ attributed to the dfDNLS  (\ref{eq:DNLS}),(\ref{eq:icsDNLS}), has  a unique global attractor ${\cal A}_{\gamma}\subset B_{r}\subset l^2$.}\label{Theorem:attractorDNLS}
\end{theorem}
\vspace{0.5cm}

\section{Congruence of the attractors ${\cal A}_{\mu}$ and ${\cal A}_{\gamma}$}

At last we establish the congruence of the attractors ${\cal A}_{\mu}$ and ${\cal A}_{\gamma}$ where we assume the following:\\
(I) $\psi_0=\phi_0 \in B\subset  B_{R_{\mu}}$, \\
(II) hypothesis (\ref{eq:Cphi}) of Theorem  \ref{Theorem:closeness} holds  with $K_{\phi}\epsilon \le R_{\mu}$,\\
(III) the conditions (\ref{eq:conditionglobal}) and (\ref{eq:ALbound}) are satisfied.

Notice that (I)-(II) confines not only $\psi(t)=S_{\mu}(t)\psi_0$, but also $\phi(t)=S_{\gamma}(t)\phi_0$  to $B_{R_{\mu}}$ for all $t\ge 0$. 

\vspace*{0.5cm}
\begin{theorem}
{\it \,\,Let assumptions (I)-(III) above hold. Then the  attractors
 ${\cal A}_{\mu}$ and ${\cal A}_{\gamma}$ coincide according to
 \begin{equation}
  {\rm dist}\left({\cal A}_{\mu},{\cal A}_{\gamma}\right)=0. 
 \end{equation}\label{Theorem:congruence}}
\end{theorem}

\noindent{\bf Proof:} 
For any bounded subset $B\subset l^2$, it holds that 
\begin{eqnarray}
 {\rm dist}({\cal A}_{\mu},{\cal A}_{\gamma})&\le& {\rm dist}({\cal A}_{\mu},S_{\mu}(t)B)+   {\rm dist}(S_{\mu}(t)B,S_{\gamma}(t)B)+ {\rm dist}(S_{\gamma}(t)B,{\cal A}_{\gamma}).\nonumber
\end{eqnarray}
  As ${\cal A}_{\mu}$ attracts any bounded set $B\subseteq B_{R_{\mu}}\subset l^2$, for  
    any $\xi >0$, there is some $T_{\mu}(\xi)>0$ such that 
  \begin{equation}
  {\rm dist}\left({\cal A}_{\mu},S_{\mu}(t)B\right)=\sup_{a\in 
 {\cal A}_{\mu}}\,\inf_{\psi_0 \in B} {\rm dist}\left(a,S_{\mu}(t)\psi_0\right)_{l^2} <\frac{\xi}{3},\,\,\,\forall t\ge T_{\mu}. \label{eq:ineq1}
 \end{equation}
  Analogously, as   ${\cal A}_{\gamma}$ attracts any bounded set $B\subseteq B_{R_{\mu}}\subset l^2$ (${\cal A}_{\gamma}$ actually attracts any bounded set in $l^2$ anyway),  
  for
    any $\xi >0$, there is some $T_{\gamma}(\xi)>0$ such that 
  \begin{equation}
  {\rm dist}\left(S_{\gamma}(t)B,{\cal A}_{\gamma}\right)={\rm dist}\left({\cal A}_{\gamma},S_{\gamma}(t)B\right)
  =\sup_{a\in 
 {\cal A}_{\gamma}}\,\inf_{\phi_0 \in B} {\rm dist}\left(a,S_{\gamma}(t)\phi_0\right)_{l^2}
  <\frac{\xi}{3},\,\,\,\forall t\ge T_{\gamma}. \label{eq:ineq2}
 \end{equation}
 Let $\overline{T}=\max\{T_{\mu}(\xi),T_{\gamma}(\xi)\}$ and consider the time interval $[0,t_0]$ with $t_0\ge \overline{T}$.
In light of Theorem \ref{Theorem:closeness} we have that for  every $t_0>0$, there exists a small $\epsilon_0>0$ and some $C>0$ such that
for every $\epsilon \in (0,\epsilon_0)$ and for all 
$\psi_0=\phi_0 \in B \subseteq B_{R_{\mu}}$, fulfilling hypothesis (II), it holds that for every $t\in [0,t_0]$
 \begin{equation}
{\rm dist}(S_{\mu}(t)B,S_{\gamma}(t)B)=\inf_{\psi_0 \in B}  {\rm dist}_{l^2}\left(S_{\mu}(t)\psi_0,S_{\gamma}(t)\psi_0\right)
  <C\cdot {\epsilon^3}.\label{eq:ineq3}
 \end{equation}

 Combining (\ref{eq:ineq1}),(\ref{eq:ineq2}) and (\ref{eq:ineq3})
we get 
\begin{eqnarray}
  0\le {\rm dist}\left({\cal A}_{\mu},{\cal A}_{\gamma}\right)&\le& 
  \frac{\xi}{3}+C\epsilon^3+\frac{\xi}{3}=\frac{2\xi}{3}+C\epsilon^3,\,\,\,\forall \xi>0,\,\,\,\,\,\,\forall \epsilon\in (0,\epsilon_0).\nonumber
\end{eqnarray}
By setting $\epsilon_0^3=\xi/(3C)$ one gets 
\begin{equation}
 {\rm dist}\left({\cal A}_{\mu},{\cal A}_{\gamma}\right)< \xi,\nonumber
\end{equation}
from which by the arbitrariness of $\xi$ it follows that 
${\rm dist}\left({\cal A}_{\mu},{\cal A}_{\gamma}\right)=0$,
and the proof is finished.

\vspace*{0.5cm}

\hspace{16.5cm} $\square$

\section{Outlook}\label{section:outlook}
Finally as an outlook on future studies,  
regarding the analytical closeness results   
there remains the problem of obtaining estimates in our statements  that hold  uniformly for any finite time.  
Moreover, extensions of the main closeness results to higher dimensional lattices $\mathbb{Z}^N$, for $N\geq 2$ and for generalized nonlinearities are of interest. As examples, we will consider the closeness of the solutions of higher dimensional discrete nonlinear Schr\"odinger lattices with generalized power and saturable nonlinearity, to those of the $N$-dimensional generalization of the AL lattice \cite{trio}.

Another aspect is the persistence of localised wave forms, supplied by the analytical solutions of the AL equation:\\
(i) under the impact of forcing and damping in the AL equation itself,\\
(ii) and in other (nonintegrable) discrete nonlinear Schr\"odinger equations in the conservative and unforced limit as well as with the inclusion of damping and forcing. The 
corresponding closeness theorems can be formulated along the lines given in this paper. Especially with view to applications the persistence of soliton solutions in (damped and forced)  nonintegrable discrete nonlinear Schr\"odinger equations plays an important role \cite{trio}. 

Furthermore,  utilising the tools provided in this manuscript, 
the 
 asymptotic features of different discrete versions of forced and damped continuum Ginzburg-Landau (GL) equations 
 represented in combined form by  
 \begin{equation}
 \frac{d u_n}{dt}=u_n+(1+i\epsilon)
 (u_{n+1}-2u_n+u_{n-1})-(1+i\epsilon)\,|u_{n}|^2(\gamma u_n+\mu (u_{n+1}+u_{n-1})),\,\,\,n\in {\mathbb{Z}}\label{eq:gGLE}
 \end{equation}
can be explored. For $\gamma=0, \mu\ne 0$ ($\gamma \neq, \mu=0$), a discrete GL equation with nonlocal (local) nonlinear term results. 
Application of our analytical closeness and congruence methods links rigorously these discrete GL equations (\ref{eq:gGLE})    
and their associated discrete nonlinear Schr\"odinger counterparts (the dfAL and dfDNLS with $\kappa=-1$ in (\ref{eq:AL}) and (\ref{eq:DNLS}), respectively) arising in the limit $\epsilon \rightarrow 0^+$ from  (\ref{eq:gGLE}). In particular a  continuity statement can be formulated proving that the solutions of the GL equations converge to those of the DNLSs. 
Furthermore, with respect to the global attractor congruence results the limit behaviour of a global attractor of a discrete GL equation can be treated by proving its upper semicontinuity  in the {\it inviscid} limit, that is  
as $\epsilon \rightarrow 0^+$ \cite{duo}. 
In addition,  in a similar way as represented in this paper for the two  discrete nonlinear Schr\"odinger equations  the congruence of the global attractors for the nonlocal GL equation and its local GL counterpart can be demonstrated \cite{prepDirk}.

\vspace*{0.5cm}

 \centerline{{\bf Acknowledgement}}
 
 I am very grateful to Nikos I. Karachalios for many fruitful discussions.

\end{document}